# A Study on Hand Proof for The Four-Color Theorem

X.-J. Wang[1*] & T.-Q. Wang[2]

[1]*State Key Laboratory of Applied Optics, Changchun Institute of Optics & Fine Mechanics and Physics, Chinese Academy of Sciences, Changchun, Jilin, 130033, China*

[2]*Department of Mechanical and Manufacturing Engineering, Miami University, Oxford, Ohio, 04560, USA*

**For the four-color theorem that has been developed over one and half centuries, all people believe it right but without complete proof convincing all[1-3]. Former proofs are to find the basic four-colorable patterns on a planar graph to reduce a map coloring[4-6], but the unavoidable set is almost limitless and required recoloring hardly implements by hand[7-14]. Another idea belongs to formal proof limited to logical operation[15]. However, recoloring or formal proof way may block people from discovering the inherent essence of a coloring graph. Defining creation and annihilation operations, we show that four colors are sufficient to color a map and how to color it. We find what trapped vertices and boundary-vertices are, and how they decide how many colors to be required in coloring arbitrary maps. We reveal that there is the fourth color for new adding vertex differing from any three coloring vertices in creation operation. To implement a coloring map, we also demonstrate how to color an arbitrary map by iteratively using creation and annihilation operations. We hope our hand proof is beneficial to understand the mechanisms of the four-color theorem.**

In 1852, Guthrie put forward a conjecture that one could always color a map in at most four colors to differently color adjacent countries with common borderlines[1-3]. The Primary proof of this conjecture focused on whether or not it was right. Kempe proved the conjecture was true in 1879 by developing the unavoidable set and the recoloring method[1,4], but Heawood turned out that Kempe's proof was not complete[5,6]. Until 1976, Appel and Haken proved the four-color theorem by using computers, improving unavoidable sets and finding a new discharging procedure[7-11]. After the development of Appel-Haken proof, scientists implemented a simplifying and understandable proof[12-14]. In 2005, Gonthier announced a formal proof for the four-color problem by using the general-purpose theorem-proving software (called Coq)[15]. However, it has not been discovered what the essence of the four-color map problem is and how to prove it completely by hand, since we must answer whether or not any PG is four-colorable and how to color a map in the wanted adjacent list.

Here, the discussed graphs are finite, undirected, simple ones without loops and multiple edges. The issue for a coloring graph is that at most how many colors, called chromatic numbers, are necessarily used.

Usually, the number of edges connected to the vertex is defined as its degree. A face is defined as an area surrounded by at least three edges. A planar graph (PG) is that any edges of a graph do not cross each other [1-3]. A defined maximal planar graph (MPG) or saturated PG is where any new edge is not able to be added on, or in which any three vertices form a triangulated face[1-2].

Only three PGs are complete (CPG) in Fig. 1, where CPG means any vertices connected. Also, it has been proved that there is no CPG graph with five vertices[1-3]. $K_3$ is the minimum triangulated MPG, but it is few used. So, we discuss graphs having at least four vertices unless otherwise noted. Besides, to implement a map mapping to PG one to one, we confine the PG is connected, which means any vertex having at least a path to arrive at any other.

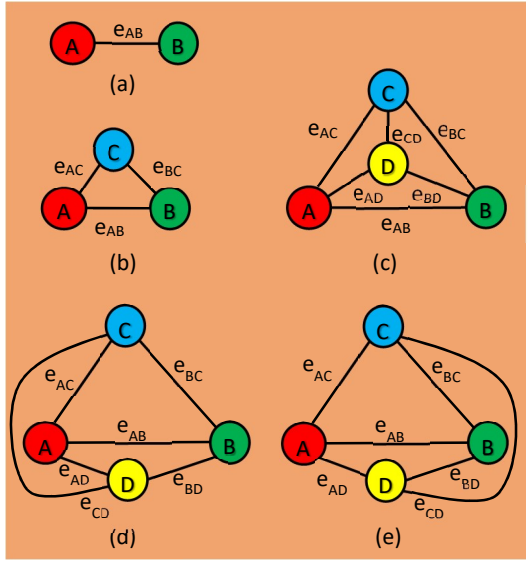

**Figure 1** Basic planar graphs. (a), CPG $K_2$; (b), CPG $K_3$; (c), CPG $K_4$; (d) and (e), isomorphic CPG $K_4$

We define, (1), the boundaries of a PG are the edges not encompassed by the finite closed area covering all vertices of the graph, which is a borderline shared by the closed area and the infinite face outside this closed area, such as $e_{AB}$, $e_{AC}$, or $e_{BC}$ in Fig. 1(c); (2), the boundary-vertex is located on the boundaries, such as vertex A, B, or C in Fig. 1(c); (3), the trapped vertex that is encompassed by all boundaries, such as vertex D that is surrounded by boundary $e_{AB}$, $e_{AC}$, and $e_{BC}$ in Fig. 1(c). Besides, the boundary connects two boundary-vertices.

Boundaries or Trapped and boundary-vertices are relative depending on how to configure boundaries. If $K_4$ in Fig. 1(c) is redrawn as its isomorphic PGs, the boundaries are $e_{CD}$, $e_{BD}$, and $e_{BC}$ and the trapped vertex becomes to vertex A in Fig. 1(d); the boundaries are $e_{CD}$, $e_{AD}$, and $e_{AC}$ and the trapped vertex is vertex B in Fig. 1(e).

$K_4$, an MPG(4) or a saturated PG, is four-coloring in Fig. from 1(c) to 1(e).

In generalization, we can find at least a subgraph $K_4$, $K_3$, or $K_2$ in the adjacency matrix or adjacency list of an MPG; otherwise, keeping its vertices unchanged for any PG, we utmost add probable edges to form an MPG.

To find the chromatic number and coloring order for an MPG, we define inside creation operation as inserting one vertex inside a triangulated face on an MPG and adding three edges from it to three vertices of this face; similarly, define outside creation operation as inserting one vertex in the infinite face on an MPG and adding three edges from it to the boundary-vertices.

For the inside creation operation in Fig. 2(a), we insert a new vertex E into a triangulated ACD of Fig 1(c) and connect vertex E to three vertices with edge $e_{EA}$, $e_{EC}$, and $e_{ED}$ in Fig. 2(a). So, three new triangulated faces replace the former face ACD keeping triangulated all faces in Fig. 2(a). The induced PG is an MPG.

Obviously in Fig. 2(a), inserting vertex E inside triangulated ACD, three boundaries and three boundary-vertices of the MPG keep unchanged, but the trapped vertices turn into vertex D and E in Fig. 2(a).

In this case, vertex E can be colored in the fourth color differing from coloring vertex A, C, or D, which vertex E is the same color as vertex B in Fig. 2(a). The induced MPG in Fig. 2(a) is four-colorable. Besides, the inserting vertex E inside the triangulated ABD or BCD of Fig. 1(c) is similar to that inside face ACD, but the results are the same.

**Lemma 1:** *There are three boundaries and three boundary-vertices on the induced MPG by using inside creation operation.*

**Proof 1:** *Since three MPG boundaries and three boundary-vertices on $K_4$ keep unchanged, there are three boundaries and three boundary-vertices on an induced MPG no matter how many vertices are added by the inside creation operations. So, Lemma 1 is workable.*

**Lemma 2:** *Induced MPG by using inside creation operation is four-colorable.*

**Proof 2:** *For an inside creation operation, it is always available to color a new adding vertex in the fourth color differing from three coloring vertices of the triangulated face where it is* inserted. So, Lemma 2 is workable.

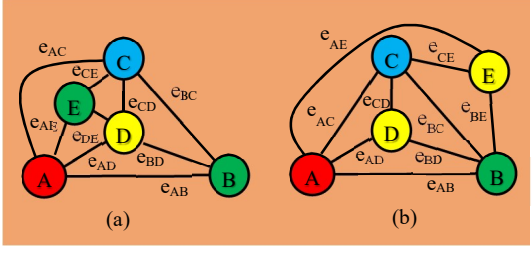

**Figure 2** Either inside or outside creation operation to insert a vertex in PG $K_4$. (a), inside; (b), outside.

For outside creation operation in Fig. 2(b), we insert a new vertex E outside the boundary of the MPG but beside boundary $e_{BC}$ of Fig. 1(c) and connect vertex E to three boundary-vertex A, B, and C with edge $e_{EA}$, $e_{EB}$, and $e_{EC}$ in Fig. 2(b). Besides, two new triangulated faces increase to this PG keeping triangulated all faces. The induced PG is an MPG.

In this case, adding vertex E becomes a boundary-vertex, vertex C turns into the trapped vertex in Fig. 2(b). The boundary-vertices are vertex E, A, and B, which number is three, unchanged.

The former boundaries $e_{AC}$ and $e_{BC}$ connecting to new trapped vertex C become to be inner edges, but, the two adding edges, $e_{AE}$ and $e_{BE}$ connecting vertex E to two former boundary-vertex A and B, turn into two new boundaries $e_{AE}$ and $e_{BE}$ in Fig. 2(b). Besides, the former boundary $e_{AB}$ connecting the two boundary-vertex A and B is still a boundary of induced MPG in Fig. 2(b). So, there are still three boundaries in the induced MPG. The number of boundaries is three, unchanged.

In this situation, vertex E can be colored in the fourth color differing from coloring vertex A, B, and C, in which vertex E is the same color as vertex D in Fig. 2(b). The induced MPG in Fig. 2(b) is four-colorable. The inserting vertex E beside boundary $e_{AC}$ or $e_{AB}$ of Fig. 1(c) is similar to that beside boundary $e_{BC}$, but the results are the same.

**Lemma 3:** *There are three boundaries and three boundary-vertices on the induced MPG by using outside creation operation.*

**Proof 3:** *By using the outside creation operation, inserting vertex E beside $e_{BC}$ for $K_4$ shown in Fig. 2(b), three boundaries are $e_{AE}$, $e_{BE}$ and $e_{AB}$, and three boundary-vertices are A, B, and E. Lemma 3 holds. Next, if inserting a vertex V beside $e_{BE}$ for Fig. 2(b), three boundaries are $e_{AV}$, $e_{BV}$ and $e_{AB}$, and three boundary-vertices are A, B, and V. Lemma 3 holds also.*

*Since the trapped vertices do not function to outside creation operations, this induced structure is like $K_4$-boundaries with three boundary-vertices.*

*Supposing until the n-th operation, Lemma 3 holds. We do the (n+1)-th operation on the n-th $K_4$-boundaries structure. Then, there are not only two new boundaries and one of one n-th boundary, but also, one new boundary-vertex and two of the n-th boundary-vertices. Lemma 3 holds also.*

**Lemma 4:** *Induced MPG by using outside creation operation is four-colorable.*

**Proof 4:** *The induction is similar to that used in Lemma 2. The coloring vertex in every outside creation operation can be in the fourth color differing from three coloring boundary-vertices.*

**Lemma 5:** *There are three boundaries and three boundary-vertices on induced MPG by using either inside or outside creation operation.*

**Proof 5:** *Combing Lemma 1 and 3, we may deduce Lemma 5.*

**Corollary 5.1:** *Any finite MPG consists of trapped vertices, inner edges, three boundaries, and three boundary-vertices.*

**Proof 5.1:** *Supposing the MPG is induced by using a series either inside or outside creation operations starting at $K_4$, the MPG poses three boundaries and three boundary-vertices according to Lemma 5. Also, the other vertices are trapped vertices except for boundary-vertices; the other edges are inner edges except for boundaries.*

**Lemma 6:** *Induced MPG by either inside or outside creation operation can be colored in four colors.*

**Proof 6:** *Combing Lemma 2 and 4, we may deduce Lemma 6.*

**Corollary 6.1:** *There is at least a way to color a finite MPG in four colors meeting the requirement of adjacent vertices in distinct colors.*

**Proof 6.1:** *Combing lemma 2 and 4, we may color the vertex in each creation operation in the fourth color.*

Corollary 6.1 does not imply that any MPG may be in four colors unless it is available to neglect the adjacent list. However, it shows having a way to color an MPG in four colors. Even if induced MPGs may be finite, we are not able to give a positive answer that unlimited various MPGs whether or not coloring in the wanted adjacent list, which challenges all proofs of the four-color theorem.

**Lemma 7:** *There is at least one three-degree vertex which degree is three for induced MPG.*

**Proof:** *Supposing the MPG is induced by using a series either inside or outside creation operations starting at $K_4$, the degree of the last adding vertex is three, which is either trapped or boundary vertex.*

**Corollary 7.1:** *Any finite MPG has a three-degree vertex which is either trapped or boundary vertex.*

**Proof 7.1:** *A finite MPG may be induced by adding every vertex one by one from the vertex set by using either inside or outside creation operations when ignoring the relationship in the adjacent list. According to Lemma 7, Corollary 7.1 holds.*

The next object is how to decrease one vertex keeping an MPG, which is inverse to creation operations.

We define an annihilation operation as only removing the three-degree vertex and three edges connecting to this vertex, which differs from removing vertex operation in graph theory; define an inside annihilation operation as removing vertex is a trapped one; define an outside annihilation operation as removing vertex is a boundary one.

For a three-degree trapped vertex on an MPG, when we remove this vertex and three edges, three triangulated faces are removed but the three vertices closest to it form a new triangulated face on induced MPG. Two triangulated faces are removed, referring to Fig. 2(a).

For a three-degree boundary-vertex on an MPG, when we remove this vertex and two boundaries and one inner edge, two triangulated faces are removed, referring to Fig. 2(b).

By using Euler's Formula[1,2], the following equality holds for any planar graph,

$$F + V = E + 2 \quad (1)$$

Where F, V, and E are the number of faces, vertices, and edges of PG respectively.

Supposing for V = n, have

$$F_n + n = E_n + 2 \quad (2)$$

For V = n-1, after one annihilation operation

$$F_{n-1} + n - 1 = E_{n-1} + 2 \quad (3)$$

Three edges are annihilated.

$$E_{n-1} = E_n - 3 \quad (4)$$

By using (1)-(4),

$$F_{n-1} + n - 1 = E_n - 3 + 2 \quad (5)$$

$$F_{n-1} + 2 + n = E_n + 2 \quad (6)$$

Then,

$$F_n = F_{n-1} + 2 \quad (7)$$

Two triangulated faces are removed after one annihilation operation, so Lemma 8 is sound.

**Lemma 8:** *After an annihilation operation, induced PG keeps MPG properties unchanged.*

**Proof 8:** *refer to (1)-(7).*

By the way, according to Lemma 7, since induced PG is an MPG after an annihilation operation, so induced MPG has three boundaries and three boundary-vertices.

We note that the annihilation operation is used only in MPG. Thus, we can repeatedly apply it to decrease vertices on an MPG until $K_4$ appears.

**Lemma 9:** *Annihilation operation ensures $MPG(V_i)$ to $MPG(V_i-1)$ until induced MPG encounters $K_4$.*

**Proof 9:** *According to Lemma 7, 8, and Corollary 7.1, we can find at least a three-degree vertex on the MPGs after annihilation operations one by one, and vertices on the MPGs are decreased one by one. Since annihilation operations inversely correspond to creation operations one to one, when the induced MPG is the same as in Fig. 2(a) or 2(b), the next induced MPG is $K_4$. The procedure refers to Fig. 4(a) and 4(b).*

Besides, Dirac proved that there was at least a $K_4$ by using his critical graph model for a four-coloring graph[16]. We know induced MPG is a four-coloring graph, Lemma 9 is true by using Dirac's proof.

To map a plane map to a graph in one-to-one correspondence, the properties of mapping relationships are described as follows[1-3]: (1), *on a map, each country is one closed area surrounded by borderlines, on a graph, each country corresponds to a vertex; (2), on a map, adjacent countries pose a common borderline on a graph, the edge denotes the borderline of adjacent countries on the map; (3), how many neighbors encompassing one country on a map equals how many edges are connected to the vertex corresponding to this country on the graph.*

Based on the above mapping properties, supposing there is a two-edge crossing on the graph which corresponds to the two-borderline crossing of two pairs of adjacent countries. On the map, two pairs of adjacent countries, four countries, have a common point. So, this common point is mapped to one edge that corresponds to only one of the two borderlines. Then, there is no edge crossing in the mapping graph. Thus, the mapping graph is a PG.

However, how to color a map is still an issue due to how to keep the adjacent list, even if we know how to color the MPG in four colors according to Corollary 7.1.

Obviously, for any one edge on a coloring MPG removed, the two coloring vertices are unnecessary to be changed or recolored in the same, because removing one edge does not destroy the other connected vertices in different colors, so long as keeping properties of the connected graph.

So, we utilize mapping $PG(V_i)$ to keep following the adjacent list of the map and form the corresponding $MPG(V_i)$ by adding edges. According to Lemma 7, 8, and 9, we can find at least one three-degree vertex on iterative $MPG(i)$.

In general, if unable to find the trapped 3-degree-vertex, we remove three boundaries to release this MPG to the PG and reselect possible three boundaries to find the 3-degree-vertex which is either trapped or new boundary-vertex. Based on Lemma 7 and Corollary 7.1, this re-selection keeps the properties of an MPG unchanged.

The best selection is to keep at least one three-degree trapped vertex in the procedure of adding edges by reselecting three different boundary-vertices in all possible ones.

Then, we iteratively apply either inside or outside annihilation operations to the three-degree vertices on $MPG(i)$, until the iterative $MPG(i-1)$ is $MPG(4)$, $K_4$, and record the order for each annihilation operation.

Inversely, we color $K_4$ in four colors and flip the order of annihilation operations from end to end as the processing order for iterative creation operations. In every step of creation, we color a new adding vertex in the fourth color differing from the other three coloring vertices until obtaining the MPG as the same size as the mapping PG, then remove adding edges and color the map in one-to-one correspondence with the mapping PG.

In a short, we exert the preliminary

verification and validation as examples.

It is well-known that $K_5$ or $K_{3,3}$ is not a PG, but $k_5$-e or $k_{3,3}$-e is. Especially, $k_5$-e is MPG, where labeling '-e' means removing any one edge from $k_5$ or $k_{3,3}$.

Some $k_5$-e are plotted in Fig. 2(a) and 2(b), either by inside or outside creation operation to $K_4$.

According to Lemma 6, $k_5$-e is four-colorable. The map corresponding to MPG $k_5$-e is four-colorable.

Another example in discussing bipartite graph, especially, $k_{3,3}$-e corresponds to a generalized geographic map including some undivided areas such as oceans, etc., although $k_{3,3}$, well-known, maybe two-colorable. So, we show two different ways how to color a map in the one-to-one correspondence from the four-coloring PG, but not aiming at using as few colors as possible.

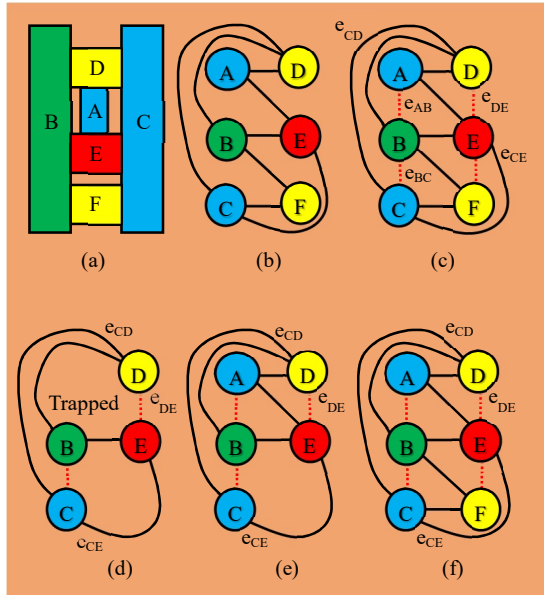

**Figure 3** Schematic for four-coloring $k_{3,3}$-e. (a), the map; (b), mapping to $k_{3,3}$-e; (c), forming MPG; (d), subgraph $K_4$ stopping annihilation; (e), creation to color vertex A; (f), creation to color vertex F.

The direct method to color $k_{3,3}$-e is heuristic. We map Fig 4(a) to Fig 4(b). For graph $k_{3,3}$-e in Fig. 3(b), add the red dashed edges to form the MPG in Fig. 3(c). After annihilations operations for inner vertex A and F, obtain subgraph $K_4$ in Fig 4(d).

According to Lemma 6, the first inside creation is to insert vertex A inside triangulated face BCD in Fig 4(d). Vertex A in Fig. 3(e) is colored in blue differing from vertex colors of triangulated BCD. The second step is to insert vertex F in Fig. 3(e). Vertex F can be colored in yellow differing from vertex colors of triangulated BCE in Fig. 3(f).

After removing added dashed edges, the coloring PG is in Fig. 3(b). Then, the coloring map one to one corresponds to this PG in Fig. 3(a).

In a generalized way, we record adding edges to form the mapping PG to the MPG in Fig. 3(e). We use annihilation to remove inner three-degree vertex F, obtain Fig. 4(a). Next, we remove the three-degree boundary-vertex C in Fig. 4(a), obtain Fig. 4(b) that is k4, then stop the operations. We record the annihilation order in the form as [in(F)_BCE, out(C)_BDE].

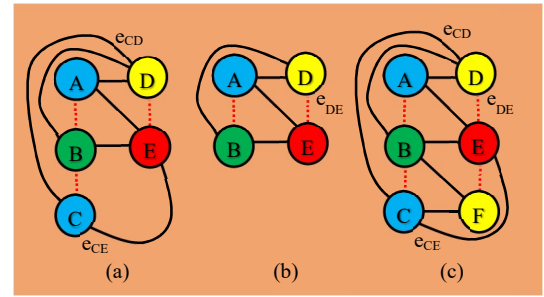

**Figure 4** Generalizing way to find $K_4$ by annihilation operations. (a), inside annihilating inner vertex F; (b), outside annihilating boundary-vertex C; (c), four-coloring forming MPG.

Inverse to creation operations, we inversely flip the annihilation order as [out(C)_BDE, in(F)_BCE] to color vertices from Fig. 4(b). We apply an outside creation to insert vertex C in triangulated BDE and color C in blue differing from colors of vertex B, D, and E in Fig. 3(e). Next, we apply an inside creation to insert vertex F in triangulated BCE and color F in yellow differing from colors of vertex B, C, and E in Fig. 4(c). So, after removing red dashed edges, coloring PG is in

Fig. 3(b). Then, we obtain the same coloring map as in Fig. 3(a).

Especially, after removing the adding edges from the four-coloring MPG, the coloring PG may degenerate to two colors, such as vertex B changes in blue, E in yellow in Fig.3 (b). So, the chromatic number of a map is at most four and the four-color theorem is true.

Our manuscript demonstrates that the essence of the four-color map problem is that there are three boundaries and the three boundary-vertices in any MPG, which decides how many colors to color a map.

Based on finding at least one three-degree vertex as an onset for any MPG, we iteratively apply annihilation operations until encountering $K_4$, and adopt its inverse order of annihilations to color the map by iteratively using creation operations.

This means that we do not need to find Kempe chains or other patterns that computers require, even for formal proof. Inversely, computers may be run by using a simple program in a short time to color a map in four colors.

No words could thoroughly review or praise a historic conjecture until it is proved by hand proof, since there may be still invisible God's hand controlling the essence[17].

**Author's contribution** All authors equally contributed to design, proving, and writing the manuscript.

**Acknowledgments** X.W. thanks Foundation of SKLAO, Innovation Foundation of CAS, and the NSF of China. T.W. acknowledges support from the fellowship of MME, Miami University, Oxford.

**Competing interests statement** The authors declare that they have no competing financial interests.

**Correspondence** and requests for materials should be addressed to X.W. (xjwang@ciomp.ac.cn) or T.W (wangt5@miamioh.edu).